\documentclass[a4paper]{amsart}
\usepackage{amssymb, amsmath, mathrsfs}
\usepackage{tikz}
\usepackage{mathtools}
\usepackage{xcolor}
\mathtoolsset{showonlyrefs}

\usepackage{todonotes}
\usepackage[
 pdfauthor={Olga Balkanova and Dimitry Frolenkov and Morten S. Risager},
 pdftitle={Prime geodesics and averages of the Zagier $L$-series},
 pdfproducer={LaTeX2e with hyperref},
 pdfcreator={Pdflatex},
 pdfdisplaydoctitle =true]
 {hyperref}
\usepackage{pdfsync}

\newtheorem{thm}{Theorem}[section]
\newtheorem{prop}[thm]{Proposition}
\newtheorem{lem}[thm]{Lemma}

\newtheorem{lemma}[thm]{Lemma}
\newtheorem*{rem*}{Remark}

\theoremstyle{definition}

\newtheorem{rem}{Remark}

\def \zL {\mathscr{L}}
\def \R {{\mathbb R}}
\def \G {\Gamma}

\def \N {\mathbf{N}}
\def \psl  {{\hbox{PSL}_2( {\mathbb R})} }
\def \pslz  {{\hbox{PSL}_2( {\mathbb Z})} }

\def \GmodH {{\Gamma\backslash\mathbb H}}
\def \e {{\varepsilon}}

\newcommand{\abs}[1]{\left\lvert #1 \right\rvert}

\providecommand{\sym}{\operatorname{sym}}
\DeclareMathOperator{\acosh}{acosh}

\title{Prime geodesics and averages of the Zagier $L$-series}

\author{Olga Balkanova}
\address{Steklov Mathematical Institute of Russian Academy of Sciences, 8 Gubkina st., Moscow, 119991, Russia}
\email{balkanova@mi-ras.ru}

\author{Dmitry Frolenkov}
\address{Steklov Mathematical Institute of Russian Academy of Sciences, 8 Gubkina st., Moscow, 119991, Russia}
\email{frolenkov@mi-ras.ru}

\author{Morten S. Risager}
\address{Department of Mathematical
  Sciences, University of Copenhagen, Universitetsparken 5, 2100
  Copenhagen \O , Denmark}
\email{risager@math.ku.dk}

\date{\today}

\begin{document}
\begin{abstract} The Zagier $L$-series encode data of real quadratic fields. We study the average size of these $L$-series, and prove asymp\-to\-tic expansions and omega results for the expansion. We then show how the error term in the asymptotic expansion can be used to obtain error terms in the prime geodesic theorem.
\end{abstract}

\maketitle
\section{Introduction}
Let $\G$ be a cofinite discrete subgroup of $\psl$ and consider the function
\begin{equation}
\Psi_\Gamma(X)=\sum_{N(P)\leq X}\Lambda(P),
\end{equation}
where the sum is over closed geodesics $P$ on  $\GmodH$ with norm $N(P)\leq X$, and $\Lambda(P)=\log(N(P_0))$ where $P$ is a power of the primitive geodesic $P_0$.
The prime geodesic theorem states that for a cofinite discrete subgroup $\G$ of $\psl$ we have
\begin{equation}\label{best-known}
\Psi_\Gamma(X)=X+O(X^{\delta+\e})
\end{equation} for some $\delta<1$.

There are several non-trivial bounds for $\delta$: Bounds for general groups e.g. \cite{Huber:1960}, \cite[Thm 10.5]{Iwaniec:2002a}, and stronger bounds for congruence groups \cite{Iwaniec:1984a}, \cite{LuoSarnak:1995}, \cite{Cai:2002} with the current record due to  Soundararajan and Young \cite{SoundararajanYoung:2013} who proved \begin{equation}\label{sy-bound}
\delta=2/3+\theta/6,\end{equation} where $\theta$ is a subconvexity exponent for quadratic $L$-functions in the conductor-aspect. See also \cite{BalkanovaFrolenkov:2019}. Concerning the value $\theta$ Conrey and Iwaniec has proved that $\theta=1/6$ is admissible \cite[Cor. 1.5]{ConreyIwaniec:2000}, and the Lindel\"of hypothesis predicts that $\theta=0$ is admissible.

For the rest of the paper we let $\G=\pslz$. One reason for the interest in the error term in this case, is the  striking relation between $\Psi_{\G}$ and averages of class number of real quadratic fields ordered by the size of the regulator, see \cite{Sarnak:1982a}.

Both \cite{SoundararajanYoung:2013} and \cite{BalkanovaFrolenkov:2019} use the relation to real quadratic fields. This relation is encoded in the Zagier $L$-series
\begin{equation}\zL_n(s)= \frac{\zeta(2s)}{\zeta(s)}\sum_{q=1}^\infty \frac{1}{q}\sum_{\substack{1\leq r\leq 2q\\r^2\equiv n\, (4q)}}1, \quad \textrm{ when }\Re(s)>1.\end{equation}
This series admits meromorphic continuation to the entire complex plane. If $n$ is non-zero and not a full square it is an entire function. We refer to \cite[\S 4]{Zagier:1977}, \cite[Sec. 4]{BalkanovaFrolenkov:2018e} for some of its basic properties.

There is a beautiful identity (see \cite[Prop 2.2]{SoundararajanYoung:2013}) linking $\Psi_\Gamma$ and the Zagier $L$-series stating that
\begin{equation}
  \Psi_\Gamma(X)=2\sum_{3\leq n\leq X}\sqrt{n^2-4}\zL_{n^2-4}(1).
\end{equation}
A relation of this kind was first discovered by Kuznetsov \cite[Eq. 7.2]{Kuznetsov:1978a}. In this paper we investigate further the relation between these two objects.

In \cite{BalkanovaFrolenkov:2017b} Balkanova and Frolenkov studied smooth averages of $\zL_{n^2-4}(s)$. Their work suggests that $\zL_{n^2-4}(1/2+it)$ has density function
\begin{equation*}m_t(x):=\begin{cases}
  \frac{1}{2\zeta(3/2)}\left(\log(x^2-4)-\frac{\pi}{2}+3\gamma -2\frac{\zeta'(3/2)}{\zeta(3/2)}-\log 8\pi\right) & \text{ if }t=0\\
\frac{\zeta(1+2it)}{\zeta(3/2+it)}+\frac{2^{1/2+it}\sin(\pi (1/2+it)/2)}{\pi^{it}}\frac{\zeta(it)}{\zeta(3/2-it))}\Gamma(it)(x^2-4)^{-it}& \text{ if }t\neq0.
\end{cases}\end{equation*}
Indeed a modification of their work (see Section \ref{sec:asymptotics}) shows the following result:
\begin{thm}\label{asymptotic}The following asymptotic formula holds
  \begin{equation}\label{compromise}
    \sum_{2<n\leq X}\zL_{n^2-4}(1/2+it)=\int_{2}^{X}m_t(u)du+O(X^{\alpha+\e})
  \end{equation}
uniformly for $\abs{t}\leq X^\e$. Here $\alpha=2(1+\theta)/3$, where $\theta$ is a subconvexity exponent for quadratic $L$-functions in the conductor-aspect.
\end{thm}
In \cite[p.6]{BalkanovaFrolenkov:2017b} it was noted that if, for $\abs{t}\leq X^\e$, the numbers $\zL_{n^2-4}(1/2+it)$ are bounded on average in windows of the form $X<n\leq X+T$ for $T\gg X^{2/3}$, then the error-term estimate in the prime geodesics theorem \eqref{best-known} holds with $\delta=2/3$. In this paper we go further and investigate what happens if we have good asymptotics in  \eqref{compromise}:

\begin{thm}\label{conditional-bound} Assume that \eqref{compromise} holds for some $\alpha >0$. Then  $\delta=1/2+\alpha/4$ is valid i.e. \begin{equation}
\Psi_\Gamma(x)=x+O(x^{1/2+\alpha/4+\e}).
\end{equation}
\end{thm}
Using the error term from Theorem \ref{asymptotic} recovers the best known bound \eqref{best-known} in the prime geodesics theorem. Considering if this can be improved it is tempting to speculate what the best possible error term is in \eqref{compromise}. We show below that if  $\alpha=1/2$ is admissible then it is optimal. Hence the limit of what might be achieved using Theorem \ref{conditional-bound} is $$\delta=5/8.$$ This is far from the conjectured $\delta=1/2$, but spectacularly better than what we know even on the generalized Lindel\"of hypothesis, which would only give $\delta=2/3$. We note that the exponent 5/8 has been proven to hold in the mean average by Cherubini and Guerreiro \cite[Thm. 1.4]{CherubiniGuerreiro:2018}. Very recently this was improved to 7/12 by Balog, Bir{\'o}, Harcos, and  Maga \cite[Thm. 1]{BalogBiroHarcosMaga:2019}.

To see that  $\alpha=1/2$ would indeed be optimal we prove the following $\Omega$-result:
\begin{thm}\label{omega}We have
  \begin{equation}
    \sum_{2<n\leq X}\zL_{n^2-4}(1/2)=\int_{2}^{X}m_0(u)du+\Omega(X^{1/2})
  \end{equation}
as $X\to\infty$.
\end{thm}

As mentioned above it is expected that
\begin{equation}\label{conjecture1}
  \Psi_\Gamma(x)=x+O(x^{1/2+\e}).
\end{equation}
This would indeed be optimal (See \cite[Theorem 3.8, p. 477]{Hejhal:1983a}).

Petridis and Risager proposed the following conjecture supporting the correctness of \eqref{conjecture1}: For any $\e>0$ and $X\geq 2$ we have
\begin{equation}\label{conjecture2}
  \sum_{0<t_j\leq T} X^{it_j}\ll T(TX)^\e,
\end{equation}
where $\lambda_j=1/4+t_j^2$ are the eigenvalues of $\Delta$. Here $\Delta$ denotes the automorphic hyperbolic Laplacian for $\pslz$. Iwaniec \cite[Lemma 1]{Iwaniec:1984a} showed that for $1\le T\le X^{1/2}\log^{-2}X$ the error term in the prime geodesic theorem can be expressed as follows
\begin{equation}\label{PrimeGeodesic to spec.sum}
\Psi_\Gamma(X)-X=2X^{1/2}\Re\left(\sum_{0<t_j\le T}\frac{X^{it_j}}{1/2+it_j}\right)+
O\left(\frac{X}{T}\log^2x\right).
\end{equation}
 Using summation by parts on the sum we easily see that the conjecture \eqref{conjecture2} implies \eqref{conjecture1}.

The following additively twisted Selberg--Linnik conjecture was suggested by Iwaniec \cite[p.139]{Iwaniec:1984a},  \cite[p.189]{Iwaniec:1984}: For any $\e>0$ and $C, D, n \geq 1$ we have \begin{equation}\label{conjecture3}
\sum_{c\leq C}\frac{e(D/c)}{c}S(n,n,c)\ll (nCD)^\e.
\end{equation} We show that this conjecture is stronger than \eqref{conjecture2}:
\begin{prop}\label{iw>>pr}
  Iwaniec's conjecture \eqref{conjecture3} implies the conjecture \eqref{conjecture2} on the spectral exponential sum.
\end{prop}

\begin{rem} The paper is organized as follows. In Section \ref{sec:asymptotics} we prove Theorem \ref{asymptotic} extending the method from \cite{BalkanovaFrolenkov:2017b}. Then in Section \ref{sec:conditional} we combine techniques from \cite{Iwaniec:1984a}, \cite{LuoSarnak:1995} and \cite{BalkanovaFrolenkov:2019a}  to prove Theorem \ref{conditional-bound}. In Section \ref{sec:omega-results} we prove the omega result in Theorem \ref{omega}. We do so by first showing that it suffices to show omega results for a smoothed problem, and then using an expansion from \cite{BalkanovaFrolenkov:2017b} combined with a lemma due to Ivi\'c and Motohashi \cite{IvicMotohashi:1990a}, to prove the desired result. Maybe surprisingly, in order to make the argument work we need the existence of a symmetric square $L$-function with non-vanishing central value. This is ensured by an asymptotic formula  for the first moment of symmetric square $L$-functions, see  \cite[Theorem 7.1.1]{Ng:2016a}, \cite{Tang:2012} and \cite{Balkanova:2019}.
\end{rem}

\section{Asymptotics for the average over $\zL_{n^2-4}(1/2+it)$}\label{sec:asymptotics}
In this section we prove Theorem  \ref{asymptotic}. For $t=0$ the result was proved in \cite[Theorem 1.3]{BalkanovaFrolenkov:2017b}, the only difference being that the proof of \cite[Theorem 1.3]{BalkanovaFrolenkov:2017b} identified the asymptotics of the integral in the main term.

The case $t\neq 0$ is similar. In particular, the main term can be evaluated using the proof of \cite[Theorem 6.3]{BalkanovaFrolenkov:2017b}. In order to estimate the error term, we follow the proof of \cite[Theorem 7.4]{BalkanovaFrolenkov:2017b}. We use the same notation as in \cite[Theorem 7.4]{BalkanovaFrolenkov:2017b}. To this end,  we define $\omega(x)$ to be a smooth characteristic function of the interval $(X,2X)$ with the length of smoothing being equal to $T$. We assume also that $\omega^{(k)}(x)\ll_k T^{-k}$. Note that we assume that $X^{\e_0}<T<X^{1-\e_0}$ since the final choice of $T$ will be $T=X^{2/3-4\theta/3}$.

Our first step is to estimate the function $h_1(\omega;1/2+it;r)$ defined as follows (see \cite[Lemma 5.2]{BalkanovaFrolenkov:2017b}):
\begin{multline}\label{eq:hwsr}
h_1(\omega;s;r)=\int_{0}^{\infty}\frac{\omega(x)}{2x^s}\frac{\sin{\pi(s/2-1/2+ir)}}{\sin{(\pi i r)}}\Gamma(1/2+s/2+ir)\\ \times \frac{\Gamma(s/2+ir)}{\Gamma(1+2ir)}\left( \frac{x}{2}\right)^{-2ir}
F(s/2+ir,s/2+1/2+ir,1+2ir;4/x^2)dx.
\end{multline}
We need to estimate $h_1(\omega;s;r)$ for both positive and negative $r$, but using that  $ h_1(\omega; s, r)=\overline {h_1(\omega; \overline s, -r)}$ we may assume that $r>0$.

Assume first that $r\ll X^{3\e}.$ In this case we approximate the hypergeometric function in \eqref{eq:hwsr} by 1 using its series representation
\begin{equation*}
F(s/2+ir,s/2+1/2+ir,1+2ir;4/x^2)=1+O\left(\frac{1+\abs{t/2+r}^2}{(1+\abs{r})x^2}\right),
\end{equation*}
for $x$ sufficiently large (See e.g. \cite[p. 38]{Good:1981b}).
Estimating all other factors in \eqref{eq:hwsr} by absolute value, we obtain
\begin{equation}\label{hestimate1}
h_1(\omega;s;r)\ll  \frac{\sqrt{X}}{\abs{r}^{1/2}}.
\end{equation}

Now assume that \begin{equation}\label{secondrange}r\gg X^{3\e}>t^3.\end{equation} In this case we can use the integral representation \cite[15.6.1]{OlverLozierBoisvertClark:2010} for the hypergeometric function in order to apply the saddle point method. Take $\alpha:=t/r\ll X^{-2\e}$ and $Z:=x^2/4$. Using \cite[15.6.1]{OlverLozierBoisvertClark:2010} we show that
\begin{multline*}
h_1(\omega;1/2+2it;r)=
\frac{\sin{\pi(-1/4+i(r+t))}}{2^{3/2+2it}\sin{(\pi i r)}}\frac{\Gamma(1/4+i(r+t))}{\Gamma(1/4+i(r-t))}\times\\
\int_{0}^{\infty}\omega(x)\int_{0}^{1}
\frac{y^{-1/4+i(r+t)}(1-y)^{-3/4+i(r-t)}}{(Z-y)^{1/4+i(r+t)}}dydx.
\end{multline*}

Consequently, bounding the factor in front of the integral by Stirling's asymptotics we find
\begin{equation}\label{hest2}
h_1(\omega;1/2+2it;r)\ll
\left|\int_{0}^{\infty}\omega(x)\int_{0}^{1}
\frac{y^{-1/4}(1-y)^{-3/4}}{(Z-y)^{1/4}}\exp(irf(\alpha,Z,y))dydx\right|,
\end{equation}
where
\begin{equation}\label{fdef}
f(\alpha,Z,y):=(1+\alpha)\log y+(1-\alpha)\log(1-y)-(1+\alpha)\log(Z-y).
\end{equation}
For simplicity, let us denote $f(y):=f(\alpha,Z,y)$ and the integral over $y$ in \eqref{hest2} by $I(\alpha,Z)$. For $X$ sufficiently large there is only one saddle point of the function $f(y)$ (the solution of  the equation $f'(y)=0$ ) in $[0,1]$ and it is located at
\begin{equation}\label{fsadle}
y_0=\frac{Z-(Z^2-Z(1-\alpha^2))^{1/2}}{1-\alpha}=\frac{Z(1+\alpha)}{Z+(Z^2-Z(1-\alpha^2))^{1/2}}.
\end{equation}
We remark that since $Z$ is  large and $\alpha$ is small, the saddle point is located near the point $1/2$. To localize the saddle point, we introduce the smooth partition of unity:
\begin{equation*}
\beta_1(y)+\beta_2(y)+\beta_3(y)=1
\end{equation*}
such that $\beta_1(y)=0$ if $y>2\delta$, $\beta_3(y)=0$ if $y<1-2\delta$, $\beta_2(y)\neq 0$ if $\delta<y<1-\delta$ for some small fixed constant $\delta$. Accordingly,
\begin{equation*}
I(\alpha,Z)=\sum_{j=1}^{3}I_j(\alpha,Z),
\end{equation*}
where
\begin{equation*}
I_{j}(\alpha,Z)=\int_{0}^{1}\frac{y^{-1/4}(1-y)^{-3/4}}{(Z-y)^{1/4}}\exp(irf(\alpha,Z,y))\beta_j(y)dy.
\end{equation*}
Since the saddle point does not belong to the support of $\beta_1(y)$ and $\beta_3(y)$, integrating by parts $I_{1}(\alpha,Z)$ and $I_{3}(\alpha,Z)$ $n$ times, we show that
\begin{equation}I_{1}(\alpha,Z)+I_{3}(\alpha,Z)\ll_n \frac{1}{Z^{1/4}r^n}.
\end{equation}
To deal with $I_{2}(\alpha,Z)$, we apply the saddle point method in the form of \cite[Theorem 1.2]{McKeeSunYe:2017} (See also \cite[Prop 8.2]{BlomerKhanYoung:2013}, getting that
\begin{multline}\label{eq:I2(alpha,Z)}
I_{2}(\alpha,Z)=\frac{y_0^{-1/4}(1-y_0)^{-3/4}}{(Z-y_0)^{1/4}}\frac{\exp(irf(y_0))}{\sqrt{rf''(y_0)}}
\left(c_0+\sum_{j=1}^{n}\frac{c_j(\alpha, Z)}{r^j} \right)\\+O(Z^{-1/4}r^{-n-1}),
\end{multline}
where $c_0=\sqrt{2\pi}\exp(\pi i/4)$ and $c_j(\alpha,Z)\ll 1$ for $j\geq 1$.

We only describe how to estimate the contribution of the main term in \eqref{eq:I2(alpha,Z)}. All other terms can be handled in the same way.

After some simplifications, we obtain
\begin{equation}\label{hest3}
h_1(\omega;1/2+2it;r)\ll\frac{1}{r^{1/2}}
\left|\int_{0}^{\infty}\frac{\omega(x)}{x^{1/2}}C(x,\alpha)\exp(irf(y_0))dx\right|+\frac{1}{Z^{1/4}r^A},
\end{equation}
where
\begin{equation}
C(x,\alpha)=\frac{y_0^{-1/4}(1-y_0)^{-3/4}}{(1-y_0/Z)^{1/4}\sqrt{f''(y_0)}}.
\end{equation}
This function is bounded and smooth in  $x$.

In order to bound the integral we need to analyze the oscillatory part $\exp(irf(y_0))$. Using \eqref{fsadle} we write
\begin{equation}
y_0=\frac{Z(1-\alpha)}{d_1},\quad  1-y_0=\frac{d_2}{d_1},\quad  Z-y_0=\frac{d_2}{1-\alpha}
\end{equation}
with
\begin{align}
  d_1=&Z+(Z^2-Z(1-\alpha^2))^{1/2}\\d_2=&(Z^2-Z(1-\alpha^2))^{1/2}-\alpha Z\\\overline {d_2}=&(Z^2-Z(1-\alpha^2))^{1/2}+\alpha Z.
\end{align}
Using \eqref{fdef} we can therefore write
\begin{multline}
f(y_0)=(1+\alpha)\log(1-\alpha^2)+(1+\alpha)\log Z-2\log(d_1)-2\alpha\log(d_2).
\end{multline}
Using
\begin{equation}
d_2\overline d_2=(1-\alpha^2)Z(Z-1)
\end{equation}
we deduce
\begin{multline}\label{f(y0)}
f(y_0)=(1-\alpha)\log(1-\alpha^2)+(1-\alpha)\log Z\\-2\alpha\log(Z-1)
-2\log d_1+2\alpha\log(\overline{d_2}).
\end{multline}
In order to simplify the expression for $\exp(irf(y_0))$, we expand the summands of $f(y_0)$ in power series with respect to $r$. For example,
\begin{equation*}
r(1-\alpha)\log(1-\alpha^2)=r(-\alpha^2+\alpha^3-\ldots)=-\frac{t^2}{r}+\frac{t^3}{r^2}-\ldots
\end{equation*}
since $\alpha=t/r$. As a result,
\begin{equation}\label{f simpify1}
\exp(ir(1-\alpha)\log(1-\alpha^2))=1-\frac{it^2}{r}+\ldots
\end{equation}
By \eqref{secondrange} we can take sufficiently many terms in the series representation and we obtain a negligibly small error term. All remaining terms can be estimated in the same way and will be smaller. Therefore, it is enough to replace $\exp(ir(1-\alpha)\log(1-\alpha^2))$ by 1.

Note further that
\begin{align*}
\frac{d_1}{Z+\sqrt{Z^2-Z}}&=\frac{Z+\sqrt{Z^2-Z(1-\alpha^2)}}{Z+\sqrt{Z^2-Z}}\\ &=  \left(1+\frac{\sqrt{Z^2-Z}}{Z+\sqrt{Z^2-Z}}\left(\sqrt{1+\frac{\alpha^2Z}{Z^2-Z}}-1\right)\right)
\end{align*}
as is easily checked. Note that the right-hand-side is close to 1. Taking $e^{-ir\log(\cdot) }$ on this expression and doing a Taylor expansion on $\sqrt{1+\frac{\alpha^2Z}{Z^2-Z}}$  and then a Taylor expansion of $e^{-ir\log(\cdot) }$ on the term close to 1 we find that
\begin{align}\label{f simpify2}
\exp\left(-2ir\log d_1 \right)= & \exp\left(-2ir\log\left(\left(Z+\sqrt{Z^2-Z}\right)\right)\right) \\
& \cdot
\left(1+\frac{\sqrt{Z^2-Z}}{Z+\sqrt{Z^2-Z}}\frac{Z\alpha^2(-2ir)}{2(Z^2-Z)}+\ldots\right).
\end{align}

Similarly, we see that
\begin{equation}
\overline d_2=
\sqrt{Z^2-Z}\left(\frac{Z\alpha}{\sqrt{Z^2-Z}}+\sqrt{1+\frac{\alpha^2Z}{Z^2-Z}}\right),
\end{equation}
and doing the same type of Taylor expansion we find
\begin{equation}\label{f simpify3}
\exp\left(2ir\alpha\log \overline{d_2}\right)=\exp\left(2ir\alpha\log\sqrt{Z^2-Z}\right)
\left(1+\frac{2irZ\alpha^2}{\sqrt{Z^2-Z}}+\ldots\right).
\end{equation}
Combining \eqref{f(y0)} and \eqref{f simpify1}, \eqref{f simpify2}, \eqref{f simpify3},  we get an expression
\begin{equation}\label{exp f_0}
\exp\left(irf(y_0)\right)=
\exp\left(irF(x,\alpha)\right)
\left(1+\ldots\right),
\end{equation}
where
\begin{align}
  F(x, \alpha)&=(1\!-\!\alpha)\log Z\!-\!2\alpha\log(Z\!-\!1)\!-\!2\log(Z\!+\!\sqrt{Z^2\!-\!Z})\!+\!2\alpha\log\sqrt{Z^2\!-\!Z}\\
  &=\log Z-\alpha\log(Z-1)-2\log(Z+\sqrt{Z^2-Z})\\
  &=-(2\log(x+\sqrt{x^2-4})+\alpha\log(x^2-4)+(6+2\alpha)\log(2)).
\end{align}
Here we have used that $Z=x^2/4$. All the terms after 1 in the last expression of \eqref{exp f_0} will be smaller in the final analysis and we ignore them below.

Using the expression \eqref{exp f_0} in \eqref{hest3} we see that we need to bound
\begin{equation}\label{hest4}
\mathcal{I}:=\frac{1}{r^{1/2}}
\abs{\int_{0}^{\infty}\frac{\omega(x)}{x^{1/2}}C(x,r)\exp\left(irF(x, \alpha)
\right)dx}.
\end{equation}
Integrating \eqref{hest4} by parts $n$ times, we infer
\begin{equation}
\mathcal{I}\ll \frac{1}{\sqrt{r}}\left(\frac{\sqrt{X}}{r}\right)^n.
\end{equation}
Therefore, for $r>X^{1/2+\epsilon}$ the integral is negligible.
In order to estimate $\mathcal{I}$ for $r\leq X^{1/2+\epsilon}$, we specialize the choice of the characteristic function
\begin{equation}\label{omega def}
\omega(x)=\frac{1}{T\pi^{1/2}}\int_X^{2X}\exp\left(-\frac{(x-K)^2}{T^2}\right)dK.
\end{equation}
For an arbitrary $A>1$ and some $c>0$ we have (see \cite{IvicJutila:2003})
\begin{equation}
\omega(x)=\begin{cases}
1+O(x^{-A})& \text{ if } x\in ]X+cT\sqrt{\log X},2X-cT\sqrt{\log X}[\\
O((|x|+X)^{-A})& \text{ if } x\not\in
[X-cT\sqrt{\log X},2X+cT\sqrt{\log X}]
\end{cases}
\end{equation}
and $1+O(T^3(T+\min(|x-X|,|x-2X|))^{-3})$ otherwise.
Using \eqref{omega def} and making the change $x=K+yT$, we show that

\begin{equation}
\mathcal{I}\ll\frac{1}{r^{1/2}}\int_X^{2X}
\int_{-\infty}^{\infty}\exp(-y^2)\frac{C(K+yT,\alpha)}{(K+yT)^{1/2}}\exp\left(irF(K+yT, \alpha)
\right)dydK.
\end{equation}
Using the rapid decay of the function $\exp(-y^2)$, we truncate the integral over $y$ at the point $v_0:=\log{(rX)}$. After that we expand all functions under the integral sign in the Taylor series.
In particular, as $x\rightarrow \infty$, we have the following Taylor series:
\begin{equation}
F(x,\alpha)=-(6+2\alpha)\log{2}-\left( 2(\alpha+1)\log{x}+2\log{2}-\frac{2+4\alpha}{x^2}+\ldots\right).
\end{equation}
Thus we show that
\begin{equation}
\mathcal{I}\ll\frac{1}{r^{1/2}}\int_X^{2X}
\int_{-v_0}^{v_0}\frac{\exp(-y^2)}{K^{1/2}}\exp\left(2ir\log(K+yT)
\right)dydK.
\end{equation}
It follows from the Taylor series expansion of the logarithm that
\begin{equation}
r\log(K+yT)=r\log{K}+\frac{ryT}{K}+r\left( \frac{yT}{K}\right)^2+\ldots
\end{equation}
Since $|y|\ll \log{X}$ and $K\sim X$, we obtain for $T\ll X^{3/4-\e_1}\ll X^{1-\e}r^{-1/2}$ that
\begin{equation}
\mathcal{I}\ll\frac{1}{r^{1/2}}\int_X^{2X}\frac{1}{\sqrt{K}}\exp(2ir(1+\alpha)\log{K})
\int_{-v_0}^{v_0}\exp(-y^2)\exp(2ir(1+\alpha)yT/K)dydK.
\end{equation}
Due to the rapid decay of the function $\exp(-y^2)$, we can enlarge the integral over $y$ to $(-\infty, \infty)$ at the cost of a negligible error term.
Evaluating the resulting integral over $y$ we have
\begin{equation}\label{eq:mathcalI}
\mathcal{I}\ll\frac{1}{r^{1/2}}\int_X^{2X}\frac{1}{\sqrt{K}}\exp(2ir(1+\alpha)\log{K})
\exp(-(r(1+\alpha)T/K)^2)dK.
\end{equation}
Note that for $r>T^{-1}X\log{X}$ the integral above is negligible. Now consider the case $r\leq T^{-1}X\log{X}$. To estimate \eqref{eq:mathcalI}, we apply the following inequality. Since we have not been able to find a proper reference we also provide a short proof:
\begin{lemma}\label{stupid-lemma} Let $-\infty\leq a\leq b \leq \infty$. Let $p(x)$ be integrable on $[a,b]$, and $q(x)$ a differentiable real function satisfying that $q'(x)\neq 0$ in the support of $p(x)$. Assume further that $p(x)/q'(x)$ is differentiable with its derivative integrable on  $[a,b]$.

  If $]a,b[$ is the union of  $m$ intervals $[a_i,b_i]$ such that $\frac{p(x)}{q'(x)}$ is monotonic on each interval $[a_i,b_i]$ we have
\begin{equation}\label{eq:sp1}
\int_{a}^b p(x)\exp(iq(x))dx \leq 2 (m+1) \max_{x\in [a,b]}\abs{\frac{p(x)}{q'(x)}}.
\end{equation}
\end{lemma}
\begin{proof} Using integration by parts we see that
\begin{equation}\abs{\int_a^b p(x)\exp(iq(x))dx}\leq
  \abs{\left. \frac{p(x)}{q'(x)}\exp(iq(x))\right\vert_a^b}+\abs{\int_a^b \frac{d}{dx}\left(\frac{p(x)}{q'(x)}\right)\exp(iq(x))dx}.
\end{equation}
If $\frac{p(x)}{q'(x)}$ is monotonic on $[a_i,b_i]$ then
\begin{equation} \int_{a_i}^{b_i}\abs{\frac{d}{dx}\left(\frac{p(x)}{q'(x)}\right)}dx=\abs{\int_{a_i}^{b_i}\frac{d}{dx}\left(\frac{p(x)}{q'(x)}\right)dx}=\abs{\left.\frac{p(x)}{q'(x)}\right\vert_{a_i}^{b_i}}\leq 2  \max_{x\in [a,b]}\abs{\frac{p(x)}{q'(x)}},
  \end{equation}
from which the lemma follows.
\end{proof}

Using Lemma \ref{stupid-lemma}, we prove that $\mathcal{I}\ll X^{1/2}r^{-3/2}$.
Therefore,
\begin{equation}\label{hestimate2}
h_1(\omega;1/2+2it;r)\ll \begin{cases}
X^{1/2}r^{-3/2} &  r<T^{-1}X\log{X},\\
X^{-100}&T^{-1}X\log{X}\leq r\leq X^{1/2+\epsilon},\\
r^{-1/2}(X^{1/2}/r)^n & r>X^{1/2+\epsilon}.
\end{cases}
\end{equation}
Now we apply \eqref{hestimate1} for $r_j\ll X^{2\e}$ and \eqref{hestimate2} for $r_j\gg X^{2\e}$ to estimate the contribution of the discrete spectrum:
\begin{equation*}
Z_D(1/2+it)\ll \sqrt{X}\left(\frac{X}{T}\right)^{1/2+\e}.
\end{equation*}
See \cite[(1.6)]{BalkanovaFrolenkov:2017b} for the definition of $Z_D(s)$. The continuous part $Z_C(s)$ \cite[(1.4)]{BalkanovaFrolenkov:2017b}  can be estimated in the same way.  We are left to estimate the holomorphic part  $Z_H(1/2+it)$, see \cite[(1.5)]{BalkanovaFrolenkov:2017b} for the definition. With this goal, we study the function
\begin{multline}\label{gdef}
g(\omega;s;k)=2^{2k-1}\sin{\frac{\pi s}{2}}\int_{0}^{\infty}\omega(x)\Gamma(k+s/2-1/2)\\
\times\frac{\Gamma(k+s/2)}{\Gamma(2k)} x^{1-2k-s}F(k+s/2-1/2,k+s/2,2k;4/x^2)dx
\end{multline}
for $s=1/2+2it.$ To deal with the hypergeometric function, we apply an approach based on the Mellin-Barnes integrals \cite{Frolenkov:2015}. To this end, we first use the property \cite[15.8.1]{OlverLozierBoisvertClark:2010}, getting
\begin{align}\label{Ftransform}
F(k+s/2&-1/2,k+s/2,2k;\frac{4}{x^2})=\\
&\left(1-\frac{4}{x^2}\right)^{-k-s/2+1/2}F(k+s/2-1/2,k-s/2,2k;-\frac{4}{x^2-4}).
\end{align}
In the same way as in \cite[Theorem 2]{Frolenkov:2015}, we prove that
\begin{equation}\label{Fasympt}
F(k+s/2-1/2,k-s/2,2k;-\frac{4}{x^2-4})=1+O\left(\frac{k+t}{\sqrt{k}(x^2-4)^{1/2}}\right).
\end{equation}
Using \eqref{gdef}, \eqref{Ftransform} and \eqref{Fasympt}, we obtain
\begin{equation}\label{gest}
g(\omega;1/2+2it;k)\ll
\int_{0}^{\infty}\frac{\omega(x) 2^{2k}}{x^{2k-1/2}}\exp(\pi|t|)\frac{|\Gamma(k+it)|^2}{\Gamma(2k)} dx.
\end{equation}
Applying the results of \cite[section 2]{Frolenkov:2015} it can be shown that
\begin{equation}\label{Gamma est small k}
\exp(\pi|t|)\frac{|\Gamma(k+it)|^2}{\Gamma(2k)}\ll\frac{t^{2k-1}\exp(tg_2(k/t))}{\Gamma(2k)},\quad\hbox{for}\quad k\le t,
\end{equation}
\begin{equation}\label{Gamma est big k}
\exp(\pi|t|)\frac{|\Gamma(k+it)|^2}{\Gamma(2k)}\ll\frac{\Gamma^2(k)}{\Gamma(2k)}\exp(t(\pi-g_1(t/k))),\quad\hbox{for}\quad k>t,
\end{equation}
where
\begin{equation}\label{g1}
g_1(z):=\int_0^z\log(1+x^2)\frac{dx}{x^2}=2\arctan z-\frac{\log(1+z^2)}{z},
\end{equation}
\begin{equation}\label{g2}
g_2(z):=\int_0^z\log(1+x^2)dx=z\log(1+z^2)-2z+2\arctan z.
\end{equation}
Substituting  \eqref{Gamma est small k} in \eqref{gest} and using the Stirling formula \cite[5.11.3]{OlverLozierBoisvertClark:2010}, we have for $k\le t$
\begin{equation*}
g(\omega;1/2+2it;k)\ll\frac{\sqrt{k}}{t}
\int_{0}^{\infty}\omega(x)\sqrt{x}\exp(kh_2(k,t,x))dx,
\end{equation*}
where
\begin{equation}\label{h2def}
h_2(k,t,x):=2+2\log(t/k)-2\log x+g_2(k/t)t/k.
\end{equation}
Since $g_2(z)-zg_2'(z)=2\arctan z-2z$ we obtain that
\begin{equation*}
\frac{\partial}{\partial k}h_2(k,t,x)=-\frac{2t}{k^2}\arctan\frac{k}{t}<0.
\end{equation*}
consequently, for $k\le t$ there exists some constants $c_1,c_2>0$ such that
\begin{equation}\label{gest2}
g(\omega;1/2+2it;k)\ll\frac{\sqrt{k}}{t}
\int_{0}^{\infty}\omega(x)\sqrt{x}(c_1t/x)^kdx\ll X^{-c_2k}.
\end{equation}
The case of $k>t$ can be treated in the same way and leads to the same estimate. Finally, we show that $Z_H(1/2+it)\ll X^{1/2}$.

Combining the above estimates with \cite[Thm. 6.2]{BalkanovaFrolenkov:2017b} and the bound $\zL_{n^2-4}(1/2+it)\ll n^{2\theta+\e}X^{\e} $ we find
\begin{align}
  \sum_{X<n\leq 2X}\zL_{n^2-4}(1/2+it)&=\sum_{n=1}^\infty\omega(n)\zL_{n^2-4}(1/2+it)+O(X^{2\theta+\e }T)\\
\label{end-game}  &=\int_{X}^{2X}m_t(u)du+Z_C(1/2+it)+Z_H(1/2+it)\\
  &\quad\quad +Z_D(1/2+it)+O(X^{2\theta+\e }T)\\
&=\int_{X}^{2X}m_t(u)du+O(X^{1/2}(X/T)^{1/2+\e}+X^{2\theta+\e }T).
\end{align}
Choosing $T=X^{2/3-4\theta/3}$ and making the dyadic partition of unity we complete the proof of Theorem \ref{asymptotic}.

We notice that the dependence on the Lindel\"of hypothesis comes from the unsmoothing in the first line of \eqref{end-game}.

\section{Conditional improvements in the prime geodesic theorem}\label{sec:conditional}
In this section we prove Theorem \ref{conditional-bound}. By Iwaniec's explicit formula \eqref{PrimeGeodesic to spec.sum} and summation by parts it suffices to bound \begin{equation}\sum_{0< t_j\leq T}t_jX^{it_j}.\end{equation} It turns out to be useful to approximate this by a smooth sum
\begin{equation}\label{relevant-sum}\sum_{0< t_j}t_jX^{it_j}\exp(-t_j/T).\end{equation}
Define
\begin{equation}\label{phi def0}
\varphi(x)=\frac{\sinh^2\beta}{2\pi}x^2\exp(ix\cosh\beta)
\end{equation}
with
\begin{equation}\label{beta def}
\beta=\frac{1}{2}\log X+\frac{i}{2T}.
\end{equation}
The reason for this choice is that the integral transform of $\varphi$,  call it $\hat\varphi$,  appearing in the Kuznetsov trace formula, satisfies
\begin{equation}\label{test-function}\hat\varphi(t)=tX^{it}\exp(-t/T)(1+O(\abs{t}^{-1})
\end{equation}
for $X,T$ sufficiently large, and $t$ bounded away from zero. The function $\hat\varphi$ is therefore well-suited to study \eqref{relevant-sum}.

Let us introduce the following notation
\begin{equation}\label{c def}
c:=-i\cosh\beta=a-ib,
\end{equation}
\begin{equation}\label{a,b def}
\begin{cases}
a:=\sinh(\log\sqrt{X})\sin((2T)^{-1}),\\
b:=\cosh(\log\sqrt{X})\cos((2T)^{-1}).
\end{cases}
\end{equation}

We remark that
\begin{equation}\label{argc}
\arg{c}=-\pi/2+\gamma,\quad 0<\gamma, \quad  T^{-1}\ll\gamma \ll T^{-1}.
\end{equation}

 Let $h(x)$ is a smooth function supported in $[N,2N]$ for some $N>1$ such that
\begin{equation}\label{h conditions}
|h^{(j)}(x)|\ll N^{-j},\,\hbox{for}\,j=0,1,2,\ldots\, \int_{-\infty}^{\infty}h(x)dx=N.
\end{equation}

Iwaniec \cite[Lemma 8]{Iwaniec:1984a} showed that by using the Lindel\"of conjecture on average  (proved by Luo and Sarnak \cite[(5)]{LuoSarnak:1995}) one can prove \begin{equation}1=\frac{\pi^2}{12N}\left(\sum_{n\in \N}h(n)\abs{\nu_j(n)}^2- r(t_j,N)\right)\end{equation} with
\begin{equation}
  \sum_{t_j\leq T}\abs{r(t_j,N)}\ll T^{2}N^{1/2}\log^2(T).
\end{equation}
Here $\nu_j(n)$ is the $n$th Fourier coefficient of the eigenfunction of the Laplace eigenvalue $1/4+t_j^2$ (appropriately normalized). Using this and \eqref{test-function} we see that in order to understand \eqref{relevant-sum} it suffices to understand
\begin{equation} \sum_{n}h(n)\sum_{t_j}\hat\varphi(t_j)\abs{\nu_j(n)}^2.
  \end{equation}
  This can be analyzed using the Kuznetsov formula, which leads to the Kloosterman sum

  \begin{equation}\label{intermediate}
    \sum_{n}h(n)\sum_{q=1}^{\infty}\frac{S(n,n;q)}{q}
    \varphi\left(\frac{4\pi n}{q}\right).
  \end{equation}
We show that if we understand the asymptotics of \begin{equation}\sum_{2<n\leq X}\zL_{n^2-4}(1/2+it)\end{equation} then we can bound \eqref{intermediate} which leads to a bound on  \eqref{relevant-sum} which ultimately leads to a bound on the error term in the prime geodesic theorem. We start by quoting a result from \cite{BalkanovaFrolenkov:2019a} which shows a relation between the sum \eqref{intermediate} and $\zL_{n^2-4}(1)$:
\begin{lem}\label{lemma estimates on the first error}
For $N,X,T\gg1$ the following asymptotic formula holds
\begin{multline}\label{sum Kloosterman decomposition2}
\sum_{n}h(n)\sum_{q=1}^{\infty}\frac{S(n,n;q)}{q}
\varphi\left(\frac{4\pi n}{q}\right)=\frac{2
\int_0^\infty h(t)dt
}
{\zeta(2)}
\sum_{n=3}^{\infty}\zL_{n^2-4}(1)\Phi(n)\\+
O\left(N\log(NX)+N^{1/2}X^{1/4+\theta}T^{3/2}\left(1+\frac{T}{X^{1/2}}\right)\right),
\end{multline}
where
\begin{equation}\label{Phi(n)def}
\Phi(x)=\frac{\sinh^2\beta}{2\pi c^2}
\frac{1-x^2/(4c^2)}{(1+x
^2/(4c^2))^{2}}.
\end{equation}
\end{lem}
\begin{proof}
See \cite[Lemma 3.4]{BalkanovaFrolenkov:2019a}.
\end{proof}

\subsection{Relation between $\zL_{n^2-4}(1)$ and $\zL_{n^2-4}(1/2+it)$}

We denote the $q$th coefficient in the Dirichlet series of $\zL_n(s)$ by $\lambda_q(n)$ i.e. for $\Re(s)>1$ we have
\begin{equation}
  \zL_n(s)=\sum_{q=1}^\infty\frac{\lambda_q(n)}{q^s}.
\end{equation}

For any $n$ and some constant $A>0$ one has
\begin{equation}\label{eq:subconvexity}
\zL_n(1/2+it)\ll (1+|n|)^{\theta}(1+|t|)^{A},
\end{equation}
where we may take $\theta$ to be any subconvexity exponent for Dirichlet $L$-functions of real primitive characters in the conductor aspect, and $t$ any subconvexity exponent in the $t$ aspect. (See \cite[Sec. 4]{BalkanovaFrolenkov:2018e} for this and other properties of $\zL_n(s)$.)  Conrey and Iwaniec \cite{ConreyIwaniec:2000} proved that we can take $\theta=1/6+\e$, and Young \cite{Young:2017b} proved a hybrid bound $\theta=A=1/6+\e$.

For $V\geq1$ define the series
\begin{equation}\label{SV def}
S_V(n^2-4):=\sum_{q=1}^{\infty}\frac{\lambda_q(n^2-4)}{q}\exp(-q/V)
,
\end{equation} which is a smoothed out analogue of
\begin{equation}
\sum_{q\leq V}^{}\frac{\lambda_q(n^2-4)}{q}.
\end{equation}
It follows from the Mellin inversion formula of the exponential function and a shift of contour (See \cite[Eqs. (1.6)-(1.8), p. 723]{Bykovskiui:1994}, \cite[p. 116, line 2]{SoundararajanYoung:2013} for details) that for $V>0$ and $n\neq2$ we have
\begin{equation}\label{approximate func.eq.}
\zL_{n^2-4}(1)=
S_V(n^2-4)-\frac{1}{2\pi i}\int_{(-1/2)}\zL_{n^2-4}(1+s)V^s\Gamma(s)ds.
\end{equation}
Using this identity we find that \begin{multline}\label{application of approx fun.eq}
\sum_{n=3}^{\infty}\zL_{n^2-4}(1)\Phi(n)=
\sum_{n=3}^{\infty}\Phi(n)S_V(n^2-4)-\\
\frac{1}{2\pi i}\int_{(-1/2)}\sum_{n=3}^{\infty}\Phi(n)\zL_{n^2-4}(1+s)V^s\Gamma(s)ds.
\end{multline}

The first summand on the right was estimated in \cite[Eq. 3.75]{BalkanovaFrolenkov:2019a} where it was found that
\begin{equation}\label{application of approx fun.eq5}
\sum_{n=3}^{\infty}\Phi(n)S_V(n^2-4)\ll
V^{1/2}T^2\log^2V.
\end{equation}

We now consider the second summand on the right of \eqref{application of approx fun.eq}.
Due to the rapid decay of $\Gamma(-1/2+it)$ we can truncate the integral over $t$ to the range $|t|\leq X^{\e}$ at the cost of an error term of $O_B(V^{-1/2}X^{-B})$ for any $B>0$.
\begin{lem}\label{intermed}
 For $|t|\ll X^{\e}\ll T $ we have
\begin{equation}\label{eq:est:1/2}
\sum_{n=3}^{\infty}\zL_{n^2-4}(1/2+it)\Phi(n)\ll X^{\alpha/2+\e}T^{2}.
\end{equation}
\end{lem}
\begin{proof}
By partial summation we have
\begin{equation}\label{eq:abelsumm}
\sum_{n=3}^{\infty}\zL_{n^2-4}(1/2+it)\Phi(n)=-\int_{1}^{\infty}\Phi'(x)\sum_{3\leq n\leq x}
\zL_{n^2-4}(1/2+it)dx.
\end{equation}
Applying Theorem \ref{asymptotic} to the inner sum, \eqref{eq:abelsumm} can be written as follows
\begin{equation}\label{mt+et}
-\int_{1}^{\infty}\Phi'(x)M_t(x)dx+O\left(\int_{1}^{\infty}\abs{\Phi'(x)}x^{\alpha+\e}dx\right).
\end{equation}
where $M_t(x)=\int_{2}^{X}m_t(u)du$.
To estimate the error term we argue as in \cite[Eqs. 3.70, 3.71]{BalkanovaFrolenkov:2019a}, and find
\begin{multline}\label{est:mtprime0}
\int_{1}^{\infty}|\Phi'(x)|x^{\alpha+\e}dx\ll |c|^{\alpha+\e}\int_{0}^{\infty}\min\left( T^3,\frac{1}{|1-x|^3}\right)x^{\alpha+\e}dx\\
\ll |c|^{\alpha+\e}T^2\ll X^{\alpha/2+\e}T^{2},
\end{multline}
where $c$ is defined by equations \eqref{c def}, \eqref{a,b def}.
The next step is to estimate the contribution of the main term in \eqref{mt+et}. Here we use integration by parts to see that
\begin{equation}
-\int_{1}^{\infty}\Phi'(x)M_t(x)dx=\int_{1}^{\infty}\Phi(x)M_t'(x)dx+O(1).
\end{equation}
We then use that $M_t'(x)=m_t(x)$ and use the specific form of $m_t(x)$.
Combining this with  \cite[Eqs. 3.56, 3.57]{BalkanovaFrolenkov:2019a} we see that
\begin{align}\label{est:mtprime}
\int_{1}^{\infty}&\Phi(x)m_t(x)dx\ll \log(X)\int_0^\infty\abs{\Phi(x)}(\log(x)+1)dx\\
&\ll \log X\int_0^\infty\frac{((1-u^2)^2+4u^2\cos^2\gamma)^{1/2}}{(1-u^2)^2+4u^2\sin^2\gamma}\abs{1+\log(u/4\abs{c})}du.
\end{align}
Using \eqref{argc} we see that the fraction in the integrand is bounded by a constant times $\min (T^2, \abs{x-1}^{-2})$.
Using this estimates we see that the full contribution of \eqref{est:mtprime} is $\ll TX^\e$. The statement follows by combining this with \eqref{mt+et}, \eqref{est:mtprime0}.
\end{proof}

\begin{lem}\label{lemma estimate on the main error}
For $X^{\e}\ll T\ll X^{1/2}$  the following estimate holds
\begin{equation}\label{estimate on sum L(1)}
\sum_{n=3}^{\infty}\zL_{n^2-4}(1)\Phi(n)\ll X^{\alpha/4+\e}T^2.
\end{equation}
\end{lem}\begin{proof}
Using \eqref{application of approx fun.eq}, \eqref{application of approx fun.eq5}, the comment before Lemma \ref{intermed}, as well as Lemma \ref{intermed} we obtain
\begin{equation}
\sum_{n=3}^{\infty}\zL_{n^2-4}(1)\Phi(n)\ll V^{1/2}T^2\log^2V+\frac{X^{\alpha/2+\e}T^{2}}{V^{1/2}}.
\end{equation}
Choosing $V=X^{\alpha/2}$, the statement follows.
\end{proof}

\begin{proof}[Proof of Theorem \ref{conditional-bound}]
By using Lemma \ref{lemma estimates on the first error} and  Lemma \eqref{lemma estimate on the main error} we find that for $X^{\e}\ll T\leq X^{1/2}$ we have
\begin{align}\frac{1}{N}\sum_{n}h(n)\sum_{q=1}^{\infty}&\frac{S(n,n;q)}{q}
\varphi\left(\frac{4\pi n}{q}\right)\\ & \ll X^{\alpha/4+\e}T^2+\log(NX)+N^{-1/2}X^{1/4+\theta}T^{3/2}.
\end{align}

 Replacing \cite[Lemma 3.5]{BalkanovaFrolenkov:2019a} in the proof of \cite[Theorem 1.2]{BalkanovaFrolenkov:2019a}  by , we show that for $X^{\e}\ll T\leq X^{1/2}$
the following estimate holds
\begin{equation}
\sum_{t_j\ll T}t_jX^{it_j}\ll X^{\alpha/4+\e}T^2.
\end{equation}
 By using partial summation, we see that when $T\leq X^{1/2}$ we have
\begin{equation}
\sum_{0<t_j\ll T}\frac{X^{it_j}}{1/2+it_j}\ll X^{\alpha/4+\e}\log T.
\end{equation}
Finally, using this to evaluate \eqref{PrimeGeodesic to spec.sum} we may choose $T=X^{1/2-\e}$ and conclude that
\begin{equation}
E(X)=O(X^{1/2+\alpha/4+\e})
\end{equation}
which proves Theorem \ref{conditional-bound}.
\end{proof}

\section{Omega results for averages of $\zL_{n^2-4}(1/2)$}\label{sec:omega-results}
We recall that $f=\Omega(g)$  $x\to\infty$ if  $f=o(g)$ does not hold, and that $f=\Omega_\pm(g)$ as $x\to\infty$ if $\limsup_{x\to\infty}f/g>0$ and $\liminf_{x\to\infty}f/g<0$. When proving Omega results for the error function it is useful to smooth it out by integrating against a suitable test function. We now describe how this can be done.

We denote the indicator function of a set $A$ as $1_A(t)$.
Consider a bump function, i.e. a smooth non-negative function $\phi$ on $\R$ supported in $[-1,1]$ satisfying $\int_\R\phi(t)dt=1$. Let $\delta>0$ be some small parameter and define $$\phi_{\delta}(t):=\delta^{-1}\phi(\delta^{-1}t).$$ Then $\phi_{\delta}$ is supported in $[-\delta, \delta]$ but otherwise has the same characteristics as $\phi$.

We now define a smooth indicator function of an interval by the convolution
\begin{equation*}
  1^{sm}_{\delta,]y,2y]}(t)=\phi_\delta*1_{]y,2y]}(t)=\int_{\R}\phi_\delta(t-v)1_{]y,2y]}(v)dv.
\end{equation*}
It is straightforward to verify that for $y>\delta$
\begin{enumerate}
  \item  $1^{sm}_{\delta,]y,2y]}(t)\in [0,1]$
  \item $1^{sm}_{\delta,]y,2y]}(t)=0$ for $t$ outside  $]y-\delta,2y+\delta]$, and
\item $1^{sm}_{\delta,]y,2y]}(t)=1$ for $t\in]y+\delta,2y-\delta]$.
\end{enumerate}

We also consider another smooth non-negative function  $\psi$ supported in $]1/2,2]$ and satisfying $\int_{\R_+}\psi(v)\frac{dv}{v}=1$. Define
\begin{equation*}
\psi_{\delta}=\delta^{-1}\psi(v^{1/\delta}).
\end{equation*} Now $\psi_{\delta}$ is supported in $[2^{-\delta}, 2^{\delta}]$ but otherwise share the same characteristics as $\psi$.

Fix now small positive parameters $0<\delta_1, \delta_2\leq 1$, and define the function
\begin{equation}\label{omega_X}
  \omega_X(t)=\int_{\R_+}\psi_{\delta_1}\left(\frac{y}{X}\right)1^{sm}_{\delta_2,]y,2y]}(t)\frac{dy}{y}.
\end{equation}
\begin{tikzpicture}
  \draw[->] (-0.2,0) -- (10,0) node[right] {$t$};
  \draw[->] (0,-0.2) -- (0,1.5) node[above] {$\omega_X(t)$};
  \draw (2pt,1) -- (-2pt,1)
			node[anchor=east] {1};
  \draw (0,2pt) -- (0,-2pt)
			node[anchor=north] {0};
  \draw (2,2pt) -- (2,-2pt)
    			node[anchor=west, rotate=-35] {$\frac{X}{2^{\delta_1}}-\delta_2$};
  \draw (4,2pt) -- (4,-2pt)
                  node[anchor=west, rotate=-35] {$2^{\delta_1}X+\delta_2$};
  \draw (6,2pt) -- (6,-2pt)
                  node[anchor=west, rotate=-35] {$2^{1-\delta_1}X-\delta_2$};
  \draw (8,2pt) -- (8,-2pt)
                  node[anchor=west, rotate=-35] {$2^{1+\delta_1}X+\delta_2$};

  \draw (1,0) -- (2,0);
  \draw (2,0) .. controls (3,0) and (3,1) .. (4,1);
  \draw (4,1) -- (6,1);
  \draw (6,1) .. controls (7,1) and (7,0) .. (8,0);
  \draw (8,0) -- (9,0);
  \draw[dashed,very thin,color=gray] (0,1) -- (4,1);
  \draw[dashed,very thin,color=gray] (4,0) -- (4,1);
  \draw[dashed,very thin,color=gray] (6,0) -- (6,1);
\end{tikzpicture}
\begin{prop}\label{omega-prop}For $X>2$ the function $\omega_{X}$ satisfies that
  \begin{enumerate}
\item \label{one} $\omega_X(t)$ is smooth and $0\leq \omega_X(t)\leq 1$
\item \label{two} $\omega_X(t)=0$ for $t\notin[2^{-\delta_1}X-\delta_2, 2^{1+\delta_1}X+\delta_2]$
\item \label{three} $\omega_X(t)=1$ for $t\in[2^{\delta_1}X+\delta_2, 2^{1-\delta_1}X-\delta_2]$
\item \label{four} for $k\geq 1$ we have $\omega_X^{(k)}\ll_{k}\frac{1}{(\delta_1X)^k}$  with support in 2 intervals of length $O(\delta_1X+\delta_2)$.
  \end{enumerate}

\end{prop}
\begin{proof}
Observe that \eqref{one} is clear from the definition.

We then note that \eqref{two} follows from noticing that if $t\notin[2^{-\delta_1}X-\delta_2, 2^{1+\delta_1}X+\delta_2]$ then for $2^{-\delta_1}\leq y/X\leq 2^{\delta_1}$ which is the support of $\psi_{\delta_1}\left(\frac{y}{X}\right)$ we have $1^{sm}_{[y,2y]}(t)=0$.

To see \eqref{three} we note that for $2^{-\delta_1}\leq y/X\leq 2^{\delta_1}$ and $t\in[2^{\delta_1}X+\delta_2, 2^{1-\delta_1}X-\delta_2]$ then $1^{sm}_{[y,2y]}(t)=1$.

Finally the claim on the support of $\omega_X^{(k)}$ follows from \eqref{two} and \eqref{three}. To bound
$\omega_X^{(k)}$ we note that repeated integration by parts gives
\begin{align*}
  \omega_X^{(k)}(t)&=\int_{\R_+}\psi_{\delta_1}\left(\frac{y}{X}\right)\left.1^{sm}_{\delta_2,]y,2y]}\right.^{(k)}(t)\frac{dy}{y}\\
  &=(-1)^{k}\int_{\R_+}\left[y^{-1}\psi_{\delta_1}\left(\frac{y}{X}\right)\right]^{(k)}1^{sm}_{\delta_2,]y,2y]}(t)dy.
\end{align*}
It is straightforward to see that for $y$ in the support of the integrand we have
$$\left[y^{-1}\psi_{\delta_1}\left(\frac{y}{X}\right)\right]^{(k)}\ll \frac{1}{(\delta_1 X)^{k}y}H_{k, \delta_1}\left(\frac{y}{X}\right)$$
for a smooth function satisfying $\int_{0}^{\infty}H_{k, \delta_1}\left(\frac{y}{X}\right)\frac{dy}{y}\ll_k 1.$ This finishes the proof.
\end{proof}

We let
$$E_t(X):=\sum_{2<n\leq X}{\zL_{n^2-4}(1/2+it)}-\int_{2}^Xm_t(u)du$$
and consider
\begin{equation} \label{smooth-error}
  E_{\delta_1,\delta_2}(X)=\int_{\R_+}\psi_{\delta_1}\left(\frac{y}{X}\right)\int_\R\phi_{\delta_2}(v)\left(E_0(2y+v)-E_0(y+v)\right)dv\frac{dy}{y}.
\end{equation}
Note that the integrand is supported in $2^{-\delta_1}X\leq y \leq 2^{\delta_1}X$, $-\delta_2\leq v\leq \delta_2$.
It suffices to prove that this smoothed out error-term  is $\Omega(X^{1/2})$ as the following lemma shows.
\begin{lem}\label{from-unsmooth-to-smooth} If $E_{\delta_1,\delta_2}(X)=\Omega(X^{1/2})$ then $E_0(X)=\Omega(X^{1/2})$.
\end{lem}
\begin{proof}
Assume to the contrary that $E_0(X)/X^{1/2}\to 0$ as $X\to \infty$.  For a given $\e>0$ we may choose $X_0$ large enough such that for $X>X_0$
$$\abs{E_0(2y+v)-E_0(y+v)}\leq \e X^{1/2},$$
 whenever  $(y,v)$ is in the support of the integrand. Combining this with $\int_{\R_+}\psi_{\delta_1}\left(\frac{y}{X}\right)\frac{dy}{y}=\int_\R\phi_{\delta_2}(v)dv=1$ we may conclude that  $\abs{E_{\delta_1,\delta_2}(X)}\leq \e X^{1/2}$, which  contradicts $E_{\delta_1,\delta_2}(X)=\Omega(X^{1/2})$.
\end{proof}
We now show how $E_{\delta_1,\delta_2}(X)$ relates to $\omega_X$ defined in \eqref{omega_X}.
\begin{lem}\label{herewego}For $X\geq 3$ we have
  \begin{equation*}
    E_{\delta_1,\delta_2}(X)=\sum_{n=1}^\infty\zL_{n^2-4}(1/2)\omega_X(n)-\int_{0}^\infty m_0(u)\omega_X(u)du.
  \end{equation*}
\end{lem}
\begin{proof} Let $X\geq 3$. For $(y,v)$ in the support of the integrand of $E_{\delta_1,\delta_2}(X)$ we have
\begin{align*}
  E_0(2y+v)-E_0(y+v)&=\sum_{y+v<n\leq 2y+v}\zL_{n^2-4}(1/2)-\int_{y+v}^{2y+v}m_0(u)du\\
  =&\sum_{n=1}^{\infty}1_{]y, 2y]}(n-v)\zL_{n^2-4}(1/2)-\int_0^\infty1_{]y, 2y]}(u-v) m_0(u)du.
\end{align*}
Inserting this in \eqref{smooth-error}, changing order of integration and making a simple change of variables $r=u-v$ gives the result.
\end{proof}

We now want to prove that
  \begin{equation}E_{\delta_1,\delta_2}=\Omega_\pm(\sqrt{X}).\end{equation} This implies clearly that $E_{\delta_1,\delta_2}=\Omega(\sqrt{X}), $ and therefore, when  combined with Lemma \ref{from-unsmooth-to-smooth}, proves Theorem \ref{omega}.

  Applying \cite[Theorem 6.3]{BalkanovaFrolenkov:2017b} and Lemma \ref{herewego} we have
\begin{equation}\label{expression}
E_{\delta_1,\delta_2}(X)=Z_C(1/2)+Z_H(1/2)+Z_D(1/2),
\end{equation}
where
\begin{equation}\label{eq:zc}
Z_C(1/2)=\frac{\zeta(1/2)}{\pi}\int_{-\infty}^{\infty}
\frac{\zeta(1/2+2ir)\zeta(1/2-2ir)}{\left|\zeta(1+2ir\right)|^2}h(\omega_X;r)dr,
\end{equation}
\begin{equation}\label{eq:zh}
Z_H(1/2)=\sum_{k\geq 6}g(\omega_X;k)\sum_{j\leq \vartheta(k)}\alpha_{j,k}L(\sym^2u_{j,k}, 1/2),
\end{equation}
\begin{equation}
Z_D(1/2)=\sum_{j}\alpha_j L(\sym ^2 u_j,1/2)h(\omega_X;r_j),
\end{equation}
where, according to \cite[Lemmata 5.1 and 5.2]{BalkanovaFrolenkov:2017b},  the  functions $h(\omega_X;r)$, $g(\omega_X;k)$  are defined as follows:
\begin{equation}\label{eq:defh}
h(\omega_X;r_j)=h_1(\omega_X;r_j)+h_1(\omega_X;-r_j),
\end{equation}
\begin{multline}\label{eq:defh1}
h_1(\omega_X;r)=\frac{\sin{(\pi(-1/4+ir))}}{2\sin(\pi i r)} \frac{\Gamma(1/4+ir)\Gamma(3/4+ir)}{\Gamma(1+2ir)}\\ \times 2^{2ir}\int_{0}^{\infty}\frac{\omega_X(x)}{x^{1/2+2ir}}F(1/4+ir,3/4+ir,1+2ir,4/x^2)dx,
\end{multline}
\begin{multline}\label{g1def}
g(\omega_X;k)=2^{2k-3/2}\int_{0}^{\infty}\omega_X(x)\Gamma(k-1/4)\\
\times\frac{\Gamma(k+1/4)}{\Gamma(2k)} x^{1/2-2k}F(k-1/4,k+1/4,2k;4/x^2)dx.
\end{multline}
We refer to \cite[Section 2]{BalkanovaFrolenkov:2017b} for the definition of $\alpha_j$ and $\alpha_{j,k}$.

\begin{lem}\label{hypasymp1} There is an $x_0>2$ such that for $r \rightarrow \infty$  we have
\begin{multline}
F\left( \frac{1}{4}+ir,\frac{3}{4}+ir,1+2ir; \frac{4}{x^2}\right)=x^{2ir}\exp(-2ir\acosh{x/2})\\ \times \left( \frac{x^2}{x^2-4}\right)^{1/4}
\left(1+\frac{1}{16ir }\left( 1-\frac{x^2-2}{x\sqrt{x^2-4}}\right) \right)+O\left(\frac{1}{x^2r^2} \right),
\end{multline}
uniformly for  all $x>x_0$.
\end{lem}
\begin{proof}
See \cite[Corollary 7.3]{BalkanovaFrolenkov:2017b}.
\end{proof}

We can now use this asymptotic expansion to to estimate $Z_D(1/2)$:

\begin{lem}\label{lem:discspectr} We have
\begin{equation}
Z_D(1/2)=\Omega_{\pm}(X^{1/2}).
\end{equation}
\end{lem}
\begin{proof}
We must analyze
\begin{equation}Z_D(1/2)=\sum_{j}\alpha_j L(\sym ^2 u_j,1/2)h(\omega_X;r_j).\end{equation}

Let $R$ be a real positive parameter (it will chosen later as $X^{\e}$) and split the sum over $r_j$ into two parts: $r_j>R$ and $r_j\leq R$.

{\bf Assume first that $r_j>R$.}
In this case we apply the asymptotic formula from Lemma \ref{hypasymp1} in order to approximate the hypergeometric function appearing in $h_1(\omega_X; r_j)$. The contribution of the error term in Lemma \ref{hypasymp1} can be handled as follows. We use the Stirling asymptotics on the gamma factors,  and bound everything else trivially. This gives that the contribution to $h_{1}(\omega_X;r_j)$ can be bounded the following way:
\begin{equation}
\int_{0}^{\infty}\frac{\omega_X(x)}{x^{5/2}r_{j}^{5/2}}dx \ll \frac{1}{r_{j}^{5/2}(X2^{-\delta_1})^{3/2}}\ll \frac{1}{r_{j}^{5/2}X^{3/2}}.
\end{equation}
Consequently, this contributes to $Z_D(1/2)$ as
\begin{equation}
\sum_{r_j>R}\alpha_j |L(\sym^2 u_j,1/2)|\frac{1}{r_{j}^{5/2}X^{3/2}}\ll \frac{1}{R^{1/2-\e}X^{3/2}}.
\end{equation}
Note that in order to get the estimate above we used the Cauchy-Schwartz inequality and the upper bound for the second moment of symmetric-square $L$-functions (see \cite[Thm 3.2]{LuoSarnak:1995}), combined with well-known bounds on $\alpha_j$ (see \cite[Cor. 0.3]{HoffsteinLockhart:1994}).

Now we estimate the contribution of the main term in Lemma \ref{hypasymp1}. To this end, consider the integral
\begin{equation}
\int_{0}^{\infty}\omega_X(x)v(x,r)\exp(-2ir_j\acosh{x/2})dx,
\end{equation}
where \begin{equation}
v(x,r)=\frac{1}{(x^2-4)^{1/4}}\left(1+\frac{1}{16ir }\left( 1-\frac{x^2-2}{x\sqrt{x^2-4}}\right)\right).
\end{equation}
Using that the derivative of $\exp(-2ir_j\acosh{x/2})$ with respect to $x$ equals $$\exp(-2ir_j\acosh{x/2})(-2ir_j)(x^2-4)^{-1/2},$$ we may integrate by parts twice in the above integral and obtain
\begin{equation}
\frac{1}{(2ir_j)^2}\int_{0}^{\infty}\Biggl(\sqrt{x^2-4}\left(\sqrt{x^2-4} \omega_X(x)v(x,r)\right)' \Biggr)' \exp(-2ir_j\acosh{x/2})dx.
\end{equation}
Using Proposition \ref{omega-prop}
 and estimating the integral trivially, this can be bounded by
$\sqrt{X}r_{j}^{-2}\delta_{1}^{-1}$ where we have assumed that $X\delta_1\gg 1$. As a result, the main term in Lemma \ref{hypasymp1} contributes to $Z_D(1/2)$ as
\begin{equation}
\sum_{r_j>R}\frac{\alpha_j }{\sqrt{r_j}}|L(\sym^2 u_j,1/2)|\frac{\sqrt{X}}{r_{j}^2\delta_1}\ll \frac{\sqrt{X}}{R^{1/2-\e}\delta_1}.
\end{equation}

{\bf Next assume that $r_j\leq R$.}

In this case, we apply the following asymptotic formula
\begin{equation}\label{asymp:hygeom2}
F(1/4+ir,3/4+ir,1+2ir,4/x^2)=1+O\left( \frac{r}{x^2}\right)
\end{equation}
when $r/x^2\ll1$, see \cite[Eq (4.13)]{Good:1981b}.
Note that
\begin{equation}
\int_{0}^{\infty}O\left(\frac{r_j}{x^2}\right)\frac{\omega_X(x)}{\sqrt{x}}dx=O\left(\frac{r_j}{X^{3/2}} \right).
\end{equation}
Therefore, the contribution of the error term in the asymptotic formula \eqref{asymp:hygeom2} to $Z_D(1/2)$ can be bounded by
\begin{equation}
\sum_{r_j\leq R}\alpha_j |L(\sym^2 u_j,1/2)|\frac{\sqrt{r_j}}{X^{3/2}}\ll \frac{R^{5/2+\e}}{X^{3/2}},
\end{equation}
where we have again used Stirling's asymptotics on the Gamma factors.
We proceed to estimate the contribution of the main term in \eqref{asymp:hygeom2}. Therefore we consider
\begin{equation}
\int_{0}^{\infty}\frac{\omega_X(x)}{x^{\frac{1}{2}+2ir_j}}dx
= \int_{-\infty}^{+\infty}\phi_{\delta_2}(v)\int_{\R_{+}}\psi_{\delta_1}\left( \frac{y}{X}\right)\int_{0}^{\infty} \frac{1_{]y,2y]}(x-v)}{x^{\frac{1}{2}+2ir_j}}dx\frac{dy}{y}dv,
\end{equation} where we have inserted the definition of $\omega_X(x)$ and interchanged summation. Consider now the two inner integrals. Evaluating the integral in $x$ and making the change of variables $y \rightarrow Xy$, we have
\begin{equation}\label{eq:doubleint}
\frac{X^{\frac{1}{2}-2ir_j}}{\frac{1}{2}-2ir_j}\int_{\R_{+}}\psi_{\delta_1}(y)\left(\left(2y+\frac{v}{X} \right)^{\frac{1}{2}-2ir_j}-\left(y+\frac{v}{X}\right)^{\frac{1}{2}-2ir_j}\right)\frac{dy}{y}.
\end{equation}
Applying the asymptotic formula
\begin{equation}
\left(2y+\frac{v}{X} \right)^{\frac{1}{2}-2ir_j}-\left(y+\frac{v}{X}\right)^{\frac{1}{2}-2ir_j}=(2y)^{\frac{1}{2}-2ir_j}-y^{\frac{1}{2}-2ir_j}+O\left(\frac{r_jv}{X\sqrt{y}}\right),
\end{equation}
and using the definition of $\psi_{\delta_1}$ and $\phi_{\delta_2}$ we find that
\begin{equation}\label{eq:asympdoubleint}
\int_{0}^{\infty}\frac{\omega_X(x)}{x^{\frac{1}{2}+2ir_j}}dx=\frac{X^{\frac{1}{2}-2ir_j}}{\frac{1}{2}-2ir_j}(2^{\frac{1}{2}-2ir_j}-1)\int_{\R_{+}}\psi_{\delta_1}(y)y^{\frac{1}{2}-2ir_j}\frac{dy}{y}+O(X^{-\frac{1}{2}}).
\end{equation}
The error term in \eqref{eq:asympdoubleint} contributes to $Z_D(1/2)$ as $O(R^{3/2}/\sqrt{X})$. In order to estimate the contribution of the main term, it is required to evaluate the integral
\begin{equation}
\int_{\R_{+}}\psi_{\delta_1}(y)y^{\frac{1}{2}-2ir_j}\frac{dy}{y}=\int_{0}^{\infty}\frac{1}{\delta_1}\psi(y^{1/\delta_1})y^{\frac{1}{2}-2ir_j}\frac{dy}{y}.
\end{equation}
Making the change of variables $y=v^{\delta_1}$ and integrating by parts, this is equal to
\begin{equation}\label{eq:hw2}
-\frac{1}{\delta_1(\frac{1}{2}-2ir_j)}\int_{0}^{\infty}\psi'(v)v^{\delta_1(\frac{1}{2}-2ir_j)}dv.
\end{equation}
Noticing that $\overline{h_1(\omega_X;r)}=h_1(\omega_X;-r)$, as follows from the definition \eqref{eq:defh1} it follows from the above considerations that if $r\leq R$   \eqref{eq:asympdoubleint} and \eqref{eq:hw2} that
\begin{equation}\label{eq:defhwr}
h(\omega_X;r)=-2\Re\Biggl[ X^{\frac{1}{2}-2ir}H(\omega_X;r)\Biggr]+O(X^{-1/2}+R^{1/2}X^{-3/2}),
\end{equation}
where
\begin{multline}
H(\omega_X;r):=\frac{\sin(\pi(-1/4+ir))}{\sin(\pi ir)}2^{2ir-1}\frac{\Gamma(1/4+ir)\Gamma(3/4+ir)}{\Gamma(1+2ir)}\\ \times \frac{1}{(\frac{1}{2}-2ir)^2\delta_1}(2^{\frac{1}{2}-2ir}-1)\int_{0}^{\infty}\psi'(v)v^{\delta_1(\frac{1}{2}-2ir)}dv.
\end{multline}
 The contribution to $Z_D(1/2)$ from the error term in \eqref{eq:defhwr} is $R^2/\sqrt X+R^{5/2}/X^{3/2}$ so anticipating $R$ being much smaller than $X$ it suffices to show that
\begin{equation}\label{eq:maincontr}
-2\sum_{r_j \leq R}\alpha_j L(\sym^2 u_j,1/2) \Re\Biggl[ X^{1/2-2ir_j}H(\omega_X;r_j)\Biggr]=\Omega_{\pm}(X^{1/2}).
\end{equation}
In order to prove \eqref{eq:maincontr} we would like to apply \cite[Lemma 3]{IvicMotohashi:1990a} (See also \cite[Lemma 3]{Ivic:1999}).
With this goal, we first extend the sum over $r_j$ to all $r_j$ at the cost of the error term $O(X^{1/2}R^{-1/2+\e})$. Furthermore, it is required to show that there is a non-zero term in the sum \eqref{eq:maincontr}. Let $r_{j_{0}}$ be the smallest of all $r_j$ such that $L(\sym^2 u_j,1/2)\neq 0$. The existence of such $j_0$ is guaranteed by \cite[Thm 1.2]{Balkanova:2019}, \cite[Thm 1]{Tang:2012},  \cite[Thm 7.1.1]{Ng:2016a}. We are therefore left to show that the integral
\begin{equation}\int_{0}^{\infty}\psi'(v)v^{\delta_1(1/2-2ir_{j_{0}})}dv\neq 0.
\end{equation}

Integrating by parts, we need to prove that
\begin{equation}\int_{0}^{\infty}\psi(v)v^{\delta_1/2-1}\exp(-2ir_{j_{0}}\delta_1\log{v})dv\neq 0.
\end{equation}
For this, it is enough to show that the real part of the integral is non-zero. But the real part equals
\begin{equation}
\int_{1/2}^{2}\psi(v)v^{\delta_1/2-1}\cos(2r_{j_{0}}\delta_1\log{v})dv,
\end{equation}
which is positive for $\delta_1$ small enough (depending on $r_{j_0}$) as the integrand is continuous and positive.

Finally, we choose the parameter $R=X^{\e}$ such that all error terms are negligible:
\begin{equation}
R^{\e}\left(\frac{1}{R^{1/2}X^{3/2}}+\frac{X^{1/2}}{R^{1/2}\delta_1}+\frac{R^{5/2}}{X^{3/2}}+\frac{R^2}{\sqrt{X}}+\frac{R^{3/2}}{X^{1/2}}+\frac{X^{1/2}}{R^{1/2}}\right)\ll X^{1/2-\e}.
\end{equation}
\end{proof}

\begin{lem}\label{lemma: holom-cont-contrib}
For any $\e>0$ the following estimate holds
\begin{equation}\label{eq:contspect}
Z_H(1/2)+Z_{C}(1/2)\ll \frac{X^{1/2}}{\log{X}}.
\end{equation}
\end{lem}

\begin{proof}
The contribution of $Z_H(1/2)$ is negligible due to the rapid decay of $g(\omega_X;k)$. See \eqref{gest2}.

Let us consider $Z_C(1/2)$.
Similarly to the proof of Lemma \ref{lem:discspectr}, the contribution of $|r|>R$ is negligible, therefore it is sufficient to consider
\begin{equation}
\frac{\zeta(1/2)}{\pi}\int_{-\infty}^{\infty}
\frac{\zeta(1/2+2ir)\zeta(1/2-2ir)}{\left|\zeta(1+2ir\right)|^2}h(\omega_X;r)\theta\left(\frac{r}{R}\right)dr,
\end{equation}
where $\theta(y)$ is a smooth characteristic function of the interval $(-1,1)$.
In view of \eqref{eq:defhwr}, in order to prove \eqref{eq:contspect}, it is required to show that
\begin{equation}\label{eq:ingt1}
\int_{-\infty}^{\infty}
\frac{\zeta(1/2+2ir)\zeta(1/2-2ir)}{\left|\zeta(1+2ir\right)|^2}H(\omega_X;r)\theta\left(\frac{r}{R}\right)X^{-2ir}dr\ll \frac{1}{\log{X}}.
\end{equation}
To this end, we integrate by parts, getting
\begin{equation}\label{eq:ingt2}
\frac{1}{\log{X}}\int_{-\infty}^{\infty}\Biggl(
\frac{\zeta(1/2+2ir)\zeta(1/2-2ir)}{\left|\zeta(1+2ir\right)|^2}H(\omega_X;r)\theta\left(\frac{r}{R}\right)\Biggr)'X^{-2ir}dr.
\end{equation}
Estimating everything in the standard way, we prove \eqref{eq:ingt1}.
\end{proof}
To summarize we see from \eqref{expression} combined with  Lemma \ref{lem:discspectr} and Lemma \ref{lemma: holom-cont-contrib} that
\begin{equation}E_{\delta_1,\delta_2}=\Omega_\pm(\sqrt{X}).\end{equation} Combining this with Lemma \ref{from-unsmooth-to-smooth} we have proved Theorem \ref{omega}.

\section{Relation between conjectures}\label{sec:iw>>pr}
In this final section we prove Proposition \ref{iw>>pr}
Our goal is to show that the conjecture \eqref{conjecture3} yields the following estimate
\begin{equation}\label{eq:newestimate}
\frac{1}{N}\sum_{n}h(n)\sum_{q=1}^{\infty}\frac{S(n,n;q)}{q}
\varphi\left(\frac{4\pi n}{q}\right)\ll T^2(NX)^{\e},
\end{equation}
where $\varphi$ is as in \eqref{phi def0}. Indeed, replacing \cite[Theorem 1.1.]{BalkanovaFrolenkov:2019a} by \eqref{eq:newestimate} in \cite[Section 4]{BalkanovaFrolenkov:2019a}
we infer that
\begin{equation}
\sum_{t_j}t_jX^{it_j}\exp(-t_j/T)\ll T^2(TX)^{\e},
\end{equation}
which implies the conjecture \eqref{conjecture2}.

Now we proceed to prove \eqref{eq:newestimate}.
We can write
\begin{equation}
\exp(ix\cosh{\beta})=\exp(ixb)\exp(-xa),
\end{equation}
which transforms \eqref{eq:newestimate} into
\begin{multline}
\frac{8\pi}{N}\sinh^2{\beta} \sum_{n}h(n)\sum_{q=1}^{\infty}\frac{S(n,n;q)}{q}
\frac{n^2}{q^2}\exp\left( ib\frac{4\pi n}{q}\right)\exp\left( -a\frac{4\pi n}{q}\right).
\end{multline}
Applying the partial summation we have
\begin{multline}\label{exp:partsum}
\frac{8\pi}{N}\sinh^2{\beta} \sum_{n}h(n)n^2\int_{1}^{\infty}\Biggl( \sum_{q<y}\frac{S(n,n,q)}{q}\exp\left( ib\frac{4\pi n}{q}\right)\Biggr)\\ \times
\frac{\partial}{\partial y}\Biggl(\frac{1}{y^2} \exp\left( -a\frac{4\pi n}{y}\right)\Biggr)dy.
\end{multline}
Since $b\asymp \sqrt{X}/2$, conjecture \eqref{conjecture3} implies that the inner sum can be estimated as follows
\begin{equation}
\sum_{q<y}\frac{S(n,n,q)}{q}\exp\left( ib\frac{4\pi n}{q}\right)\ll (n^2\sqrt{X}y)^{\e}.
\end{equation}
Then \eqref{exp:partsum} can be bounded by a constant times
\begin{multline}
\frac{\sinh^2{\beta}}{N} \sum_{n}h(n)n^2\int_{1}^{\infty}(n^2\sqrt{X}y)^{\e}\Biggl[\frac{1}{y^3}\exp\left( -a\frac{4\pi n}{y}\right)+\frac{an}{y^4}\exp\left( -a\frac{4\pi n}{y}\right)\Biggr]dy
\\ \ll \frac{X}{N}\sum_{n}h(n)n^2\frac{(nX)^{\e}}{(na)^2}\ll T^2(NX)^{\e},
\end{multline}
where we have used that $a \asymp \frac{\sqrt{X}}{4T}$.
This concludes the proof of Proposition \ref{iw>>pr}.

\section*{Acknowledgement}
The research of Olga Balkanova was funded by RFBR, project number 19-31-60029. The research of Morten S. Risager was supported by the Independent Research Fund Denmark DFF-7014-00060B.

\bibliographystyle{amsplain}
\providecommand{\bysame}{\leavevmode\hbox to3em{\hrulefill}\thinspace}
\providecommand{\MR}{\relax\ifhmode\unskip\space\fi MR }
\providecommand{\MRhref}[2]{%
  \href{http://www.ams.org/mathscinet-getitem?mr=#1}{#2}
}
\providecommand{\href}[2]{#2}

\end{document}